# The study of odd graceful and odd strongly harmonious for bipartite graph ☆

XiuyingWang[a], Ying Zhou[a], Haihong Wang[a], Chunfeng Liu[b] *

([a]College of Information Science and Technology, Qingdao University of Science and Technology,
Shandong, Qingdao 266000)
([b]College of Science, Liaoning University of Technology, Liaoning,Jinzhou, 121013)

**Abstract:** In this paper, we investigate odd graceful graph , odd strongly harmonious graph, bipartite graph and their relationship. We proved following results: (1) if $G$ is odd strongly harmonious graph, then $G$ is odd graceful graph ;(2) if $G$ is bipartite odd graceful graph, then $G$ is odd strongly harmonious graph; (3) graphs $K^+$、$K^*$、$G_n$、$C_{4k}^r$ and $P_{r,2s-1}$ are odd strongly harmonious graphs.

## 1. Introduction

In 1969 Rose [1] developed the initial concept of graceful graph, and then Golomb[2] defined it in greater detail; In 1991 Gnanajothi[3] proposed odd graceful graph and assumed that each tree was odd graceful graph. In 1982, Frank Hsu D[4] introduced strongly harmonious label for graph. Then Liang and Bai[5] developed odd strongly harmonious labeling and the research of the graph labeling has raised. By now, many studies have been accomplished [5-17]. However, existed research achievements were for specified graphs, the generalized results and the relationship of odd graceful graph and odd strongly harmonious label are inadequate. In this paper, we investigate odd graceful graph, odd strongly harmonious graph, bipartite graph, and their relationship.

In this paper, all graphs we discussed are the finite, undirected, and simple. $(p, q)$ - graph $G$ is a graph with $p$ vertices and $q$ edges. Let $[m, n]$ denotes a set of continuous and non-negative integers $\{m, m+1, m+2,..., n\}$, where $m$ and $n$ are integers and $0 \leq m < n$; $[k, l]^o$ denotes the odd set $\{k, k+2, k+4,..., l\}$, where $k$ and $l$ are both odd and $0 \leq k < l$; $[s, t]^e$ denotes an even set $\{s, s+2, s+4,..., t\}$, where $s$ and $t$ are both even numbers and $0 \leq s < t$.

**Definition 1.1[3].** If $(p, q)$-graph $G$ is a simple graph and there is a injection $f:V(G) \to [0, 2q-1]$ subject to following conditions:

1) $\forall x,y \in V(G)$, $x \neq y$, there is $f(x) \neq f(y)$;

2) $\forall xy \in E(G)$, let $f(xy)=|f(x)-f(y)|$, there is $\{f(xy) | xy \in E(G)\} = [1, 2q-1]^o$.

---

☆ This work was Supported by the National Natural Science Foundation of China (61773107, 61603168)

* Corresponding Author at College of Science,Liaoning University of Technology, Liaoning,Jinzhou, 121013) E-mail: liuchunfeng968@163.com

Then $G$ is called an odd graceful graph, and f is called odd graceful label of $G$.

**Definition 1.2[5].** If $(p, q)$-graph $G$ is a simple graph and there is a injection $g:V(G) \to [0, 2q-1]$ subject to following conditions:

1) $\forall x,y \in V(G)$, $x \neq y$, there is $g(x) \neq g(y)$;

2) $\forall xy \in E(G)$, let $g(xy)=g(x)+g(y)$, there is $\{g(xy) | xy \in E(G)\} = [1, 2q-1]^o$.

Then $G$ is called an odd strongly harmonious graph, and $f$ is called an odd strongly harmonious label of $G$.

Let $h$ is a label of graph $G$, the set of vertex and edge for graph $G$ are denoted respectively by $h(V(G)) = \{h(x) | x \in V(G)\}$ and $h(E(G)) = \{h(xy) | xy \in E(G)\}$. And $h_{max}(V(G))=\max h(V(G))$, $h_{min}(V(G))=\min h(V(G))$, $h_{max}(E(G))=\max h(E(G))$, $h_{min}(E(G))=\min h(E(G))$.

**Definition 1.3.** Let $(p, q)$-graph $G$ be a bipartite graph, $(V_1,V_2)$ is the bipartite division of $V(G)$. If $G$ has an odd graceful label $f$ such that $\forall x \in V_1$ and $\forall y \in V_2$, thereis $f(x) < f(y)$, then $G$ is a bipartite odd graceful graph, and $f$ is a bipartite odd graceful label of $G$.

**Definition 1.4.** Let $(p, q)$-graph $G$ be a bipartite graph, $(V_1,V_2)$ is the bipartite division of $V(G)$, where $V_1=\{u_i | i \in [1, s]\}$, $V_2=\{v_j | j \in [1, t]\}$, $s+t=p=|V(G)|$.

(1) If there is a bipartite odd graceful label $f$ of $G$, such that $f(u_i) < f(u_{i+1})$, $i \in [1, s-1]$; $f(v_j) < f(v_{j+1})$, $j \in [1, t-1]$, and $\forall i \in [1, s]$ and $\forall j \in [1, t]$, thereis $f(u_i) < f(v_j)$, then $f$ is called the standard odd graceful label of $G$.

(2) If there is an odd strongly harmonious label $g$ of $G$, such that $g(u_i) < g(u_{i+1})$, $i \in [1, s-1]$; $g(v_j) < g(v_{j+1})$, $j \in [1, t-1]$, $g(u_1)=0$, then $g$ is called the standard odd strongly harmonious label of $G$.

By definition 1.1 and definition 1.4, there are the bipartite odd graceful graphs have the standard odd graceful labels and the odd strongly harmonious graphs have the standard odd strongly harmonious labels. The following lemmas can be obtained from the above definitions.

**Lemma 1.5** If $(p, q)$-graph $G$ is a bipartite odd graceful graph, and $f$ is a standard odd graceful label of $G$, then

(1) $f(u_1)=\min f(V(G))=0$, $f(v_t)=\max f(V(G))=2q-1$, $u_1v_t \in E(G)$;

(2) $f(v_1)-f(u_s)=1$, $v_1u_s \in E(G)$;

(3) If $G$ is a connected graph, then $f(u_i)$ ($i \in [1, s]$) is even, $f(v_j)$ ($j \in [1, t]$) is odd.

**Proof.** (1) Since $G$ is a bipartite odd graceful graph, by definition 1.1, $\exists xy \in E(G)$, such that $f(xy)=2q-1$, $f(x)=0$, $f(y)=2q-1$ (or $f(x)=2q-1$, $f(y)=0$). Let $f(x)=0$, $f(y)=2q-1$. By definition 1.4(1), $x=u_1 \in V_1$, $y=v_t \in V_2$.

(2) By definition 1.4(1), there is $f(u_1)<f(u_2)<\cdots<f(u_s)<f(v_1)<\cdots<f(v_{t-1})<f(v_t)$, so $\min\{f(v_i)-f(u_j) \mid i \in [1, s], j \in [1, t]\} = f(v_1)-f(u_s) \geq 1$. $\exists xy \in E(G)$, such that $f(xy)=1$, so $f(v_1)-f(u_s)=1$, and $v_1 u_s \in E(G)$.

(3) Since $f$ is a standard odd graceful label of $G$, there is $f(u_1)=0$, $f(v_t)=2q-1$. Since $G$ is a connected graph, and by definition 1.3 and 1.4(1), then $f(u_i)$ ($i \in [1, s]$) is even, $f(v_j)$ ($j \in [1, t]$) is odd. □

**Lemma 1.6.** If $(p, q)$-graph $G$ is a connected odd strongly harmonious graph, $g$ is a standard odd strongly harmonious label of $G$, then

(1) $g(u_1)=0$, $g(v_1)=1$, and $u_1 v_1 \in E(G)$;

(2) If $\exists u \in V(G)$, such that $g(u)=g_{max}(V(G))=2q-1$, then $d(u)=1$;

(3) $g_{max}(V(G)) \geq q$;

(4) If $\exists u \in V(G)$, such that $g(u)=g_{max}(V(G))=q$, then $\exists v \in V(G)$, there is $g(v)=q-1$, $uv \in E(G)$.

(5) If $G$ is connected, then $g(u_i)$ ($i \in [1, s]$) is even, $g(v_j)$ ($j \in [1, t]$) is odd.

**Proof.** (1) Since $G$ is an odd strongly harmonious graph, by definition 1.2 and By definition 1.4(2), there is the edge of labeled 1 in $G$, so $u_1 \in V_1$, $v_1 \in V_2$, such that $g(u_1)=0$, $g(v_1)=1$, $u_1 v_1 \in E(G)$.

(2) Since $G$ is an odd strongly harmonious graph, by definition 1.2, if $\exists u \in V(G)$, $g(u)=g_{max}(V(G))=2q-1$, then $u$ can only be adjacent to the vertex of labeled 0 in $G$, so $d(u)=1$.

(3) If $g_{max}(V(G))<q$, then $\forall x \in V(G)$, there is $g(x) \leq q-1$, hence $\forall e=uv \in E(G)$, there is $g(uv)=g(u)+g(v) <2q-2$, the result is conflict with definition 1.2.

(4) Let $u_0 \in V(G)$, $g(u_0)=g_{max}(V(G))=q$. Because $g_{max}(E(G))=2q-1$, there is $v_0 \in V(G)$ in $G$, such that $g(v_0)=q-1$, $u_0 v_0 \in E(G)$.

(5) Let $g$ is a standard odd strongly harmonious label of $G$, there is $g(u_1)=0$, $f(v_1)=1$. Since $G$ is a connected graph, and by definition 1.2 and definition 1.4(2), there is $g(u_i)$ ($i \in [1, s]$) is even, $g(v_j)$ ($j \in [1, t]$) is odd. □

**Lemma 1.7.** If graph $G$ is an odd graceful graph, then $G$ is a bipartite graph.

**Proof.** Let $f$ is an odd graceful label of $G$, define $V_1=\{f(u) \mid u \in V(G), f(u) \text{ is even}\}$, $V_2=\{f(v) \mid v \in V(G),$ $f(v)$ is odd}. Obviously $V(G)= V_1 \cup V_2$, and $V_1 \cap V_2=\emptyset$, since $G$ is odd graceful, by definition1.1, $\forall xy \in E(G)$, there is $x \in V_1$, $y \in V_2$(or $x \in V_2, y \in V_1$), otherwise, there is $x,y \in V_1$(or $x,y \in V_2$), there is $f(xy)=|f(x)-f(y)|$ is even, the result is conflict with definition1.1. So $G$ is a bipartite graph. □

**Lemma 1.8.** If graph $G$ is an odd strongly harmonious graph, then $G$ is a bipartite graph.

**Proof.** Let $g$ is an odd strongly harmonious label of $G$, and let $V_1=\{u|\ u \in V(G),\ f(u)$ is even$\}$, $V_2=\{v|\ v \in V(G), g(v)$ is odd$\}$. Obviously, $V(G)= V_1 \cup V_2$, $V_1 \cap V_2==\emptyset$, because $g$ is is an odd strongly harmonious label of $G$, by definition1.2, $\forall xy \in E(G)$, there is $x \in V_1$, $y \in V_2$(or $x \in V_2$, $y \in V_1$), otherwise, $x,y \in V_1$(or $x,y \in V_2$), there is $g(xy)= g(x)+g(y)$ is even, the result is conflict with definition1.2. So $G$ is a bipartite graph. □

By lemma **1.**7 and lemma**1.**8, we obtained that research the odd gracefulness and odd strong harmony of graph just need to research the bipartite graph.

## 2. Main results and proofs

**Theorem 2.1.** Let $(p, q)$-graph $G$ is a bipartite graph, if $G$ is an odd strongly harmonious graph, then $G$ is an odd graceful graph.

**Proof.** Let$(U,V)$ be the bipartite division of vertices set of $G$, $U=\{u_i|i \in [1, s]\}$, $V=\{v_j|j \in [1, t]\}$. Let $g$ is an odd strongly harmonious label of $G$, and meet definition1.4(2), by lemma1.6(1), there is $g(u_1)=0, g(v_1)=1, u_1v_1 \in E(G)$. By definition1.4(2), there is $g(u_i)(i \in [1, s])$ is even, $g(v_j)(j \in [1, t])$ is odd. A new label $f$ of $G$ is defined as follow：

$f(u_i)=g(u_i),\qquad i \in [1, s]$;

$f(v_j)=2q-g(v_j),\quad j \in [1, t]$.

According to the definition of $f$, there is $\forall x,y \in V(G)$, $x \neq y$, there is $f(x) \neq f(y)$, $\min f(V(G))=f(u_1)=g(u_1)=0$, $\max f(V(G))=f(v_1)=2q-g(v_1)= 2q-1$. We noticed $f(u_iv_j)=|f(u_i)-f(v_j)|=f(v_j)$ $-f(u_i)= (2q-g(v_j))-g(u_i)=2q-(g(v_j)+g(u_i))= 2q-g(v_ju_i)$ $(i \in [1, s], j \in [1, t])$. Since $\{g(xy)|xy \in E(G)\}=\{g(u_iv_j)|\ i \in [1, s], j \in [1, t]\}= [1, 2q-1]^o$, so $\{f(xy)|xy \in E(G)\}=\{f(u_iv_j)|\ i \in [1, s], j \in [1, t]\}=\{2q-g(v_ju_i)|i \in [1, s], j \in [1, t]\}= [1, 2q-1]^o$, then $f$ is an odd graceful label of $G$, $G$ is an odd graceful graph. □

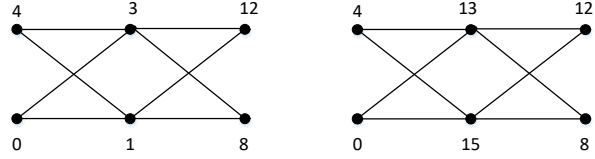

(a)The odd strongly harmonious label of graph $G_1$     (b)The odd graceful label of graph $G_1$

Fig.1 An example to explain theorem 2.1

By theorem 2.1, if $G$ is an odd strongly harmonious graph, then $G$ is an odd graceful graph. But the converse is not true. For example, the six-point circle $C_6 = v_1v_2v_3v_4v_5v_6v_1$ is an odd graceful graph($f(v_1)=0, f(v_2)=11, f(v_3)=2, f(v_4)=9, f(v_5)=4, f(v_6)=3$), but $C_6$ is not an odd strongly harmonious graph. If the conditions of odd graceful graphs are strenthened, we can obtain theorem 2.2 below.

**Theorem 2.2.** Let $(p, q)$-graph $G$ be a bipartite graph, if $G$ is a bipartite odd graceful graph, then $G$ is an odd strongly harmonious graph.

**Proof.** Let $(U,V)$ be the bipartite division of the vertices set $V(G)$ of $G$, $U=\{u_i|i\in[1, s]\}$, $V=\{v_j|j\in[1, t]\}$. Let $f$ be an odd graceful label of $G$ and meet definition1.4(1), by lemma 1.5(1), there is $f(u_1)=0, f(v_t)=2q-1, u_1v_t\in E(G)$. By Lemma 1.5(3), there is $f(u_i)(i\in[1, s])$ is even, $f(v_j)(j\in[1, t])$ is odd. A new label $g$ of $G$ is defined as follow:

$g(u_i)=f(u_i), \quad i\in[1, s]$;

$g(v_j)=2q-f(v_j), \quad j\in[1, t]$.

According to the definition of $g$, there is $\forall x,y\in V(G)$, $x\neq y$, there is $g(x)\neq g(y)$, $\min g(V(G))=g(u_1)=f(u_1)=0$, $\max g(V(G))=g(v_1)=2q-f(v_1)=2q-1$. And $g(u_iv_j)=g(u_i)+g(v_j)=f(u_i)+(2q-f(v_j))=(2q-(f(v_j))-f(u_i))=2q-f(v_ju_i)$ $(i\in[1, s], j\in[1, t])$. Because $\{g(xy)|xy\in E(G)\}=\{f(u_iv_j)| i\in[1, s], j\in[1, t]\}=[1, 2q-1]^o$, therefore $\{g(xy)|xy\in E(G)\}=\{g(u_iv_j)| i\in[1, s], j\in[1, t]\}=\{2q-f(v_ju_i)|i\in[1, s], j\in[1, t]\}=[1, 2q-1]^o$, then $g$ is an odd strongly harmonious label of $G$, $G$ is an odd strongly harmonious graph. □

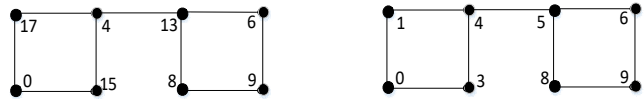

(a) The bipartite odd graceful label of graph $G_2$     (b) The odd strongly harmonious label of graph $G_2$

Fig.2 An example to explain theorem 2.2

**Theorem 2.3.** Let $(p, q)$-graph $G$ is a bipartite graph, if $G$ is an odd strongly harmonious graph, $g$ is an standard odd strongly harmonious label of $G$ and $g(u_s)+ g(v_t)< 2q$, then $G$ is a bipartite odd graceful graph.

**Proof.** Let $G$ is an odd strongly harmonious graph. $(U,V)$ be the bipartite division of the vertices $V(G)$ set of $G$, and $U=\{u_i|i\in[1, s]\}$, $V=\{v_j|j\in[1, t]\}$. Let $g$ is an standard odd strongly harmonious label of $G$ and subject to definition 1.4(2), by lemma 1.6(1), there is $g(u_1)=0$, $g(v_1)=1$, $u_1v_1\in E(G)$. By lemma 1.6(5), there is $g(u_i)(i\in[1, s])$ is even, $g(v_j)(j\in[1, t])$ is odd.

We define a new label $g$ of $G$ as follow:

$f(u_i)=g(u_i)$, $\quad i\in[1, s]$;

$f(v_j)=2q-g(v_j)$, $\quad j\in[1, t]$.

By theorem 2.1, $f$ is an odd graceful label of $G$. By the definition of $f$, there is $f(v_1)=2q-g(v_1)=2q-1$, $f(u_i)<f(u_{i+1})$ ($i\in[1, s-1]$), $f(v_{j+1})<f(v_j)$ ($j\in[1, t-1]$). By definition 1.3, we just need to to prove $f(v_t)>f(u_s)$. Because $f(v_t)-f(u_s)=(2q-g(v_t))-g(u_s)=2q-(g(v_t)+g(u_s))\geqslant 2q-((2q-1)+g(u_s))=1+g(u_s)>0$, i.e. $f(v_t)>f(u_s)$, hence $G$ is a bipartite odd graceful graph. □

**Theorem 2.4.** Complete bipartite graph $K_{m,n}$ is a bipartite odd graceful graph.

**Proof.** Let $(U,V)$ is the bipartite division for the vertices $V(K_{m,n})$ set of $K_{m,n}$, $U=\{u_i|i\in[1, m]\}$, $V=\{v_j|j\in[1, n]\}$. Obviously, there is $|E(K_{m,n})|=mn$. The label $f$ of $K_{m,n}$ is defined as follow:

$f(u_i)=2(i-1)$, $\quad i\in[1, m]$;

$f(v_j)=2mj-1$, $\quad j\in[1, n]$.

By definition of $f$, $\forall x,y\in V(K_{m,n})$, $x\neq y$, there is $f(x)\neq f(y)$, $\min f(V(K_{m,n}))=\min g(U)=f(u_1)=0$, $\max f(V(K_{m,n}))= \max f(V)=f(v_n)=2mn-1$, and $\{f(xy)|xy\in E(K_{m,n})\}= \{f(u_iv_j)|\ i\in[1, m],\ j\in[1, n]\}=\{2mj-2i+1|j\in[1, n], i\in[1, m]\}=[1, 2mn-1]^o$. So $f$ is the bipartite odd graceful label of $K_{m,n}$, and graph $K_{m,n}$ is a bipartite odd graceful graph. □

By theorem 2.2 and theorem 2.4, we can obtain following theorem 2.5.

**Theorem 2.5.** Complete bipartite graph $K_{m,n}$ is an odd strongly harmonious graph.

**Theorem 2.6.** Let $K_{m_i,n_i}(i\in[1, t])$ be $t$ complete bipartite graphs, then graph $U_{i=1}^t K_{m_i,n_i}$ (i.e. $K^+ = U_{i=1}^t K_{m_i,n_i}$) is an odd strongly harmonious graph.

**Proof.** Let $K_{m_i,n_i}(i\in[1, t])$ be $t$ complete bipartite graphs, $(U_i,V_i)$ is the bipartite division for the vertices set $V(K_{m_i,n_i})$ of $K_{m_i,n_i}$, $U_i=\{u_{i,a}|i\in[1, t], a\in[1,m_i]\}$, $V_i=\{v_{i,b}|i\in[1, t], b\in[1,n_i]\}$. Obviously, $|E(K^+)|=\sum_{k=1}^t m_k n_k$. The vertex label $g$ of $K^+$ is defined as follow:

$g(u_{i,a})=2\sum_{k=1}^{i-1} m_k+2(i-1)+2(a-1)$, $\quad i\in[1, t], a\in[1,m_i]$;

$g(v_{i,b})=2\sum_{k=1}^{i-1} m_k n_k-2\sum_{k=1}^{i-1} m_k -2(i-1)+2m_i(b-1)+1$, $\quad i\in[1, t], b\in[1,n_i]$.

According to the definition of $g$, $\forall x,y \in V(K^+)$, $x \neq y$, there is $g(x) \neq g(y)$, $\min g(V(K^+))=g(u_{1,1})=0$, $\max g(V(K^+))=g(v_{t,nt})= 2\sum_{k=1}^{t} m_k n_k - 2\sum_{k=1}^{t} m_k + 2(t-1) < 2\sum_{k=1}^{t} m_k n_k - 1 = 2|E(K^+)|-1$. By $g(u_{i,a}v_{i,b})=g(u_{i,a})+g(v_{i,b})= 2\sum_{k=1}^{i-1} m_k n_k + 2(a-1)+2m_i(b-1)+1$, $i\in[1,t]$, $a\in[1,m_i]$, $b\in[1,n_i]$, there is $\{g(xy)|xy \in E(K^+)\}=\{g(u_{i,a}v_{i,b})| i \in [1, t], a \in [1,m_i], b \in [1,n_i]\}= \{2\sum_{k=1}^{i-1} m_k n_k + 2(a-1)+2m_i(b-1)+1 \mid i \in [1, t], a \in [1,m_i], b \in [1,n_i]\}= [1, 2\sum_{k=1}^{t} m_k n_k -1]^o$, then $g$ is an odd strongly harmonious label of $K^+$, $K^+$ is an odd strongly harmonious graph. □

By theorem 2.3 and theorem 2.6, we can obtain following theorem 2.7.

**Theorem 2.7.** Let $K_{m_i,n_i}(i\in[1,t])$ is $t$ complete bipartite graphs, then their union set graph $K^+$ is a bipartite odd graceful graph.

**Theorem 2.8.** Let $K_{m_i,n_i}$ ($i\in[1,t]$) is $t$ complete bipartite graphs, $(U_i,V_i)$ is the bipartite division for the vertices set $V(K_{m_i,n_i})$ of $K_{m_i,n_i}$, $U_i=\{u_{i,j}|i\in[1,t],j\in[1,m_i]\}$, $V_i=\{v_{i,j}|i\in[1,t],j\in[1,n_i]\}$. The graph $K_{m_i,n_i}$ are merged into a new graph $K^*$ by coinciding vertices $v_{i,1}$ bonded together. Seeing Fig 3. Then graph $K^*$ is an odd strongly harmonious graph.

**Proof.** In graph $K^*$, set coinciding vertices $v_{i,1}$ as $v_1$. By the composition of graph $K^*$, there is $|E(K^+)|=\sum_{k=1}^{t} m_k n_k$. The label $g$ of $K^*$ is define as follows:

$g(v_1)=0$;

$g(u_{i,j})=2\sum_{k=1}^{i-1} m_k + 2j-1$,  $\qquad i\in[1,t], j\in[1,m_i]$;

$g(v_{1,j})=2\sum_{k=1}^{t} m_k n_k - 2m_1(j-1)$, $\qquad j\in[2,n_i]$;

$g(v_{i,j})= g(v_{i-1,2})+ g(v_{i-1,1}) - g(u_{i,1})2 - 2m_i(j-2)$, $\qquad i\in[2,t], j\in[2,n_i]$.

We can prove that $g$ is an odd strongly harmonious label of $K^*$, its proof is similar to theorem 2.6. Here is omitted.

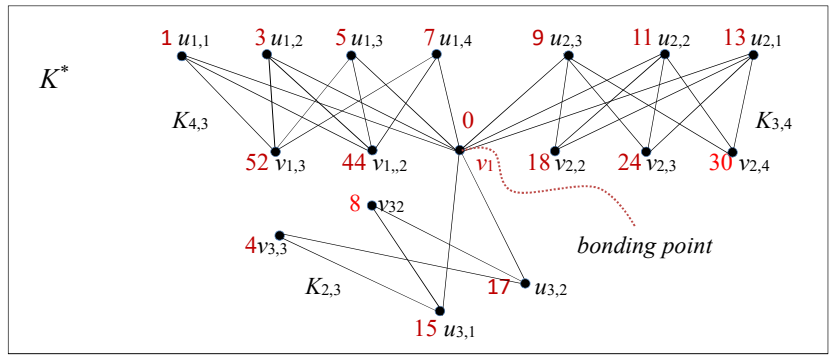

Fig.3 The graph $K^*$ made $K_{2,3}$, $K_{3,4}$ and $K_{4,3}$ by bonding point and its an odd strongly harmonious labeling

By theorem 2.3 and theorem 2.8, we can obtain following theorem 2.9

**Theorem 2.9.** Let $K_{m_i,n_i}$ ($i \in [1, t]$) is a $t$ complete bipartite graph, where the bipartite division for the set of vertices $V(K_{m_i,n_i})$ of $K_{m_i,n_i}$ is ($U_i, V_i$), $U_i = \{u_{i,j} | i \in [1, t], j \in [1, m_i]\}$, $V_i = \{v_{i,j} | i \in [1, t], j \in [1, n_i]\}$, The graph $K_{m_i,n_i}$ are merged into a new graph $K^*$ by coinciding vertices $v_{i,1}$ bonded together. Then $K^*$ is a bipartite odd graceful graph.

**Theorem 2.10.** Let $G$ is a graph, where $V(G_n) = \{u_{i,j} | i \in [1, n], j \in [1,4]\}$, $E(G_n) = \{u_{i,1}u_{i+1,1} | i \in [1, n-1]\} \cup \{u_{i,3}u_{i+1,3} | i \in [1, n-1]\} \cup \{u_{i,1}u_{i,2}, u_{i,2}u_{i,3}, u_{i,3}u_{i,4}, u_{i,4}u_{i,1} | i \in [1, n]\}$, then $G_n$ is an odd strongly harmonious graph.

**Proof.** According to composition of $G_n$, there is $|V(G_n)|=4n$, $|E(G_n)|=6n-2$. Every vertex label $g$ of $G_n$ is definded as follow:

$$g(u_{i,1}) = \begin{cases} 6i - 5, & i \equiv 1 \pmod 2, \\ 6i - 2, & i \equiv 0 \pmod 2. \end{cases}$$

$$g(u_{i,2}) = \begin{cases} 6(i-1), & i \equiv 1 \pmod 2, \\ 6i - 5, & i \equiv 0 \pmod 2. \end{cases}$$

$$g(u_{i,3}) = \begin{cases} 6i - 3, & i \equiv 1 \pmod 2, \\ 6(i-1), & i \equiv 0 \pmod 2. \end{cases}$$

$$g(u_{i,4}) = \begin{cases} 6(i-1) + 4, & i \equiv 1 \pmod 2, \\ 6i - 3, & i \equiv 0 \pmod 2. \end{cases}$$

Next we prove that $g$ is an odd strongly harmonious label of $G_n$.

(1) For the same $i \in [1, n]$, obviously there is $g(u_{i,s}) \neq g(u_{i,t})$, $1 \leq s, t \leq 4$; and max$\{g(u_{i,1}), g(u_{i,2}), g(u_{i,3}), g(u_{i,4})\} \leq$ min$\{g(u_{i+1,1}), g(u_{i+1,2}), g(u_{i+1,3}), g(u_{i+1,4})\}$, $i \in [1, n]$. So $\forall x, y \in V(G_n)$, $x \neq y$, there is $g(x) \neq g(y)$.

min$\{g(x) | x \in V(G_n)\} = g(u_{1,2}) = 0$, max$\{g(x) | x \in V(G_n)\} = g(u_{n,4}) \leq 12n-5$.

(2) Let $D_1 = \{g(u_{i,1}u_{i,2}), g(u_{i,2}u_{i,3}), g(u_{i,3}u_{i,4}), g(u_{i,4}u_{i,1}), g(u_{i,1}u_{i+1,1}), g(u_{i,3}u_{i+1,3}) | i \in [1, n-1]\}$
$= [1, 12n-13]^\circ$, $D_2 = \{g(u_{n,1}u_{n,2}), g(u_{n,2}u_{n,3}), g(u_{n,3}u_{n,4}), g(u_{n,4}u_{n,1})\} = \{12n-11, 12n-9, 12n-7, 12n-5\}$, there is $\{g(e) | e \in E(G_n)\} = D_1 \cup D_2 = [1, 12n-5]^\circ$.

By above (1) and (2), it can be verified that $G_n$ is an odd strongly harmonious graph. Seeing Fig 4.

□

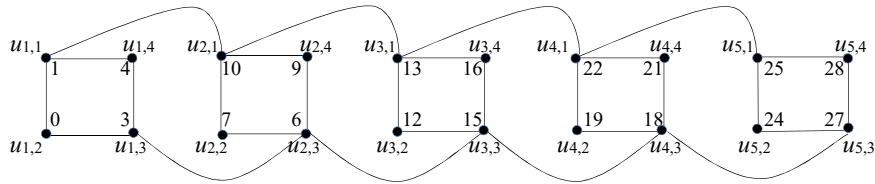

Fig.4 The graph $G_5$ and its an odd strongly harmonious labeling

**Theorem 2.11.** The $r$-crown of Circle $C_{4k}$ ($k \geq 1$) is an odd strongly harmonious graph.

**Proof.** Let $C_{4k}=u_1u_2u_3\cdots u_{4k}u_1$, $v_{i,j}$ ($1\leq i\leq 4k$, $1\leq j\leq r$). The $r$-crown of $C_{4k}$ is represented as $C_{4k}^r$ (cycle $C_{4k}$ with $r$ pendant edge attached at each vertex), i.e. $V(C_{4k}^r)=V(C_{4k})\cup\{v_{i,j}\mid 1\leq i\leq 4k, 1\leq j\leq r\}$, $E(C_{4k}^r)=E(C_{4k})\cup\{u_iv_{i,j}\mid 1\leq i\leq 4k, 1\leq j\leq r\}$, Seeing Fig.5. $|V(C_{4k}^r)|=4k(r+1)$, $|E(G)|=4k(r+1)$.

Define $g: V(G)\to[1, 2q-1]$ as follows:

$g(u_{2i-1})=2(r+1)i-2(r+1)$,     $i\in[1, k]$;

$g(u_{2i-1})=2(r+1)i-2r$;     $i\in[k+1, 2k]$;

$g(u_{2i})=2(r+1)i-1$,     $i\in[1, 2k]$;

$g(v_{2i-1,j})=2(r+1)i+2j-2(r+1)-1$,     $i\in[1, 2k], j\in[1, r]$;

$g(v_{2i,j})=2(r+1)i+2j-2(r+1)$,     $i\in[1, k], j\in[1, r]$;

$g(v_{2i,j})=2(r+1)i+2j-2r$,     $i\in[k+1, 2k], j\in[1, r]$.

By definition of $g$, it be verified easily that $g$ is an odd strongly harmonious label of graph $C_{4k}^r$, so $C_{4k}^r$ is an odd strongly harmonious graph. □

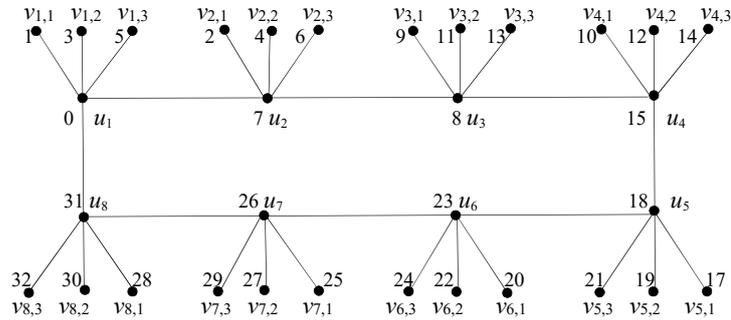

Fig.5   The graph $C_8^3$ and its an odd strongly harmonious labeling

**Theorem 2.12.** $P_{r,2s-1}$ is an odd strongly harmonious graph.($r,s$ are positive integers).

Let $x,y$ be two vertices, we defined graph $P_{a,b}$ consisting of $b$ internal disjoint paths, in which each path length is $a$, start and end vertices of path are merged with the vertices x and y respectively. For example graph $P_{7,3}$ is shown in Fig 6 below. The internal vertices of the $i$th and $r$-length path denoted by $u_0^i, u_1^i, u_2^i, \ldots, u_r^i$ in $P_{a,b}$. For all $i$, there is $u_0^i=u_0$, $u_r^i=u_r$. The $j$th vertex($j\neq 0,r$) (called as column-$j$ vertex) in every path denoted by $u_j^1, u_j^2, u_j^3, \ldots, u_j^{2s-1}$. $i=1,2,\ldots$, $j=1,2,\cdots,r-1$.

**Proof.** By the composition of $P_{r,2s-1}$, there is $|V(P_{r,2s-1})|=(r-1)(2s-1)+2$, $|E(G_n)|=r(2s-1)$. Every vertex label $g$ of $P_{r,2s-1}$ is defined as follow:

when $r\equiv 0\pmod 2$

$g(u_0)=0, \quad g(u_r)=r(2s-1),$

$g(u^i_{2j-1})= 2i+2(2s-1)j-4s+1, \quad i\in[1, 2s-1], j\in[1, \frac{r}{2}],$

$g(u^i_{2j})=4s-4i+2(2s-1)j, \quad i\in[1, 2s-1], j\in[1, \frac{r-2}{2}],$

when $r\equiv 1\pmod 2$

$g(u_0)=0, \quad g(u_r)=r(2s-1),$

$g(u^i_{2j-1})= 2i+2(2s-1)j-4s+1, \quad i\in[1, 2s-1], j\in[1, \frac{r-1}{2}],$

$g(u^{2i-1}_{2j})= 2-2i +2(2s-1)j, \quad i\in[1, s], j\in[1, \frac{r-1}{2}],$

$g(u^{2i}_{2j})= 2s -2i +2(2s-1)j, \quad i\in[1, s-1], j\in[1, \frac{r-1}{2}].$

It be verified easily that $g$ is an odd strongly harmonious label of graph $P_{r,2s-1}$, so $P_{r,2s-1}$ is an odd strongly harmonious graph. □

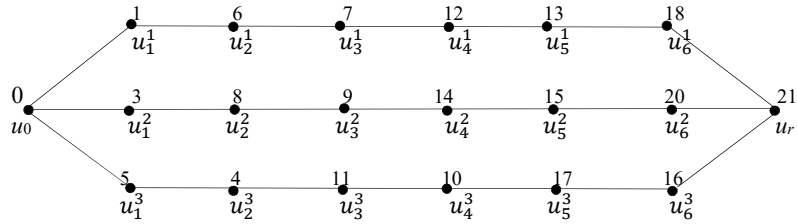

Fig.6 The graph $P_{7,3}$ and its an odd strongly harmonious labeling

By theorem 2.3, theorem 2.10 and theorem 2.12, we can obtain the following inferences:

**Inference 2.13.** Let $G$ is a graph, where $V(G_n)=\{u_{i,j}|i\in[1, n], j\in[1,4]\}$, $E(G_n)=\{u_{i,1}u_{i+1,1}|i\in[1, n-1]\} \cup \{u_{i,3}u_{i+1,3}|i\in[1, n-1]\} \cup \{u_{i,1}u_{i,2}, u_{i,2}u_{i,3}, u_{i,3}u_{i,4}, u_{i,4}u_{i,1}|i\in[1, n]\}$; then $G_n$ is a bipartite odd graceful graph.

**Inference 2.14.** The $r$-crown of circle $C_{4k}$ ($k\geq 1$) is a bipartite odd graceful graph.

**Inference 2.15**[17]. $P_{r,2s-1}$ is a bipartite odd graceful graph.($r,s$ are integers).